\documentclass[letterpaper, 10 pt, conference]{ieeeconf}  

\IEEEoverridecommandlockouts             
\let\labelindent\relax

\usepackage{graphicx,amsmath,amssymb}
\usepackage[noadjust]{cite}

\usepackage{comment,color}
\usepackage{kotex}
\usepackage{mathtools}
\usepackage{algorithm}
\usepackage{algorithmic}
\usepackage{tikz}
\usepackage{enumitem}
\usepackage{pgfplots}
\usepackage{grffile}
\usepackage{stfloats}
\pgfplotsset{compat=newest}
\usepgfplotslibrary{groupplots}   
\usetikzlibrary{plotmarks}
\usetikzlibrary{arrows.meta}
\usepgfplotslibrary{patchplots}
\usepackage{siunitx}
\usetikzlibrary{shapes,arrows,arrows.meta,positioning,calc}
\tikzstyle{block} = [draw, fill=white, rectangle, minimum height=3em, minimum width=4em]
\tikzstyle{sum} = [draw, fill=white, circle, node distance=1cm]

\usepackage{amsthm}

\newcommand{\bbR}{\ensuremath{{\mathbb R}}}

\newcommand{\bbN}{\ensuremath{{\mathbb N}}}

\newcommand{\bbE}{\ensuremath{{\mathbb E}}}

\newcommand{\bbP}{\ensuremath{{\mathbb P}}}

\newcommand{\calN}{\mathcal{N}}
\newcommand{\calP}{\mathcal{P}}

\newcommand{\calF}{\mathcal{F}}

\newcommand{\ini}{\mathsf{ini}}

\newcommand{\Tr}{\mathsf{Tr}}

\newtheorem{thm1}{\bf Theorem}
\newtheorem{prop1}{\bf Proposition}
\newtheorem{lem1}{\bf Lemma}
\newtheorem{asm1}{\bf Assumption}
\newtheorem{defn1}{\bf Definition}
\newtheorem{rem1}{\bf Remark}

\newtheorem{cor1}{\bf Corollary}

\newtheorem{prob1}{\bf Problem}

\title{\LARGE \bf
On the Invariance of Risk-Sensitive LQR Gain \\ Under Input Randomization
}

\author{Yeongjun Jang%
\thanks{*This work was supported by the National Research Foundation of Korea (NRF) grant funded by the Korea government (MSIT) (No. RS-2024-00353032 and RS-2026-25504174).
}
\thanks{Y.~Jang is with ASRI, Department of Electrical and Computer Engineering, Seoul National University, Seoul, 08826, Korea (email: jangyj0512@snu.ac.kr).
}
}

\begin{document}

\maketitle
\thispagestyle{empty}
\pagestyle{empty}

\begin{abstract}
This paper shows that the optimal gain of the risk-sensitive linear quadratic regulator (LQR) problem is invariant under input randomization when the cost penalizes the realized noisy input.
This may appear counterintuitive at first glance because certainty equivalence does not hold for risk-sensitive LQR and input randomization inflates the effective process noise.
Nonetheless, the gain is preserved because the input noise enters not only the system dynamics but also the cost functional, and its total effect eventually vanishes.
Consequently, the optimal gain and the Riccati recursion need not be recomputed, and the increment in the optimal cost can be readily evaluated in closed form.
This result facilitates the use of risk-sensitive LQR in applications that employ input randomization for privacy or exploration, such as dynamic watermarking for replay attack detection.

\end{abstract}

\section{Introduction}\label{sec:intro}

The linear quadratic regulator (LQR) is a standard problem in control theory, in which the goal is to design a policy that minimizes a quadratic cost functional for a linear dynamic system \cite{Kalm60,AndeMoor07}.
An appealing feature of LQR is that the optimal policy takes the form of linear feedback, whose gain admits a closed-form representation characterized by a backward Riccati recursion.
Moreover, it is certainty equivalent, meaning that the optimal gain does not depend on the statistics of the process noise \cite{Simo56}. 
Therefore, the optimal gain for the noise-free system can be applied to the noisy one without any loss of optimality.

However, as LQR only minimizes the expected value of the cost functional, it is inherently risk-neutral and offers no mechanism for penalizing low-probability high-cost trajectories. 
To address this limitation, risk-sensitive LQR \cite{Jaco73, Whit81}, also known as linear exponential quadratic regulator (LEQR), has been introduced.
By minimizing an expected value of the exponential quadratic cost functional, it assigns heavier weight to low-probability high-cost trajectories.

Notably, risk-sensitive LQR also admits a linear feedback optimal policy, whose gain can be computed via a backward Riccati recursion \cite{Jaco73, Whit81}.
Unlike LQR, however, this recursion depends on the process noise covariance, and certainty equivalence no longer holds. 
As a result, the optimal gain must be recomputed whenever the process noise statistics change, which can be computationally burdensome and has therefore hindered a wider usage of risk-sensitive LQR.

For example, input randomization, which refers to deliberately corrupting the nominal input with additional noise before it is applied to the plant, has been employed in various applications such as dynamic watermarking for replay attack detection \cite{MoyiSino09,SatcKuma17,NahaTeix23}, differential privacy \cite{CortDull16,YazdJone22}, and path integral control \cite{Kapp05,PatiHana24}.
The aforementioned applications have been developed primarily within the standard LQR framework, as input noise inflates the covariance of the effective process noise and is expected to alter the optimal gain in the risk-sensitive setting.
To the best of our knowledge, however, this effect has not been analyzed rigorously, motivating the central question of this paper. \textit{Is the risk-sensitive LQR gain invariant under input randomization, and if not, how much does it change?}

Interestingly, we show that the optimal gain and the associated Riccati recursion of risk-sensitive LQR are invariant under input randomization when the cost penalizes the realized noisy input.
This is not a consequence of certainty equivalence, which does not hold for risk-sensitive LQR. Instead, it arises because the input noise enters both the system dynamics and the cost functional, and its total effect on the gain eventually vanishes. 
Input randomization in risk-sensitive linear quadratic Gaussian control has also been studied in \cite{LleoRung25}.
Instead of penalizing the realized input, the quadratic cost is averaged over the input noise distribution before entering the exponential function. 
We show that the optimal gain varies under their formulation or when only the nominal input is penalized.
This highlights the role of penalizing the realized input in preserving the gain. 

Indeed, the optimal cost increases, but thanks to this invariance, the increment can be readily evaluated in closed form.
Input randomization also introduces an additional feasibility condition that is coupled with that of the original risk-sensitive LQR problem. 
We show that these conditions are monotone in the risk-sensitivity parameter, i.e., they are jointly feasible whenever the parameter is chosen sufficiently small, offering a practical parameter selection guideline.


\textit{Notation}:
Let $\bbR$ and $\bbN$ denote the sets of real numbers and positive integers, respectively.
For a sequence $v_0,\ldots, v_N$ of scalars or vectors, we define $v_{0:N}:=[v_0^\top,\ldots,v_N^\top]^\top$.
We write $\calN(\mu,\Sigma)$ to denote a multivariate Gaussian distribution with mean $\mu$ and covariance $\Sigma$.
The determinant of a square matrix $A$ is denoted by $|A|$.
We use the Loewner order, i.e., for symmetric matrices $A$ and $B$, we write
$A \succ B$ if $A - B$ is positive definite. The relations $\succeq$, $\prec$, and $\preceq$ are understood accordingly.
All random variables in this letter are defined on the probability space $(\Omega,\calF,\bbP)$, where $\Omega$ is a sample space, $\calF$ is a $\sigma$-algebra, and $\calP$ is the probability measure.

\section{Review of Risk-Sensitive LQR}\label{sec:review}

Consider a time-varying discrete time plant written by
\begin{subequations}\label{eq:LEQR}
\begin{align}\label{eq:sys}
    x_{t+1} = A_tx_t + B_t u_t + w_t,
\end{align}
where $x_t\in\bbR^n$ is the state with the initial value $x_0=x^\ini \in \bbR^n$ and $u_t\in\bbR^m$ is the input.
The process noise $w_t\in\bbR^n$ is assumed to be mutually independent zero-mean Gaussian with covariance $W_t\succ 0$, i.e., $w_t\sim \calN(0,W_t)$.

Let $N\in\bbN$ denote the terminal time. 
For each $t=0,\ldots, N-1$, we define the quadratic cost accumulated along a trajectory $(x_{t:N}, u_{t:N-1})$ by 
\begin{multline}
    \label{eq:costQuad}
    L_t(x_{t:N}, u_{t:N-1}) \coloneqq x_N^\top Q_N x_N \\
    + \textstyle\sum_{k=t}^{N-1} \left(x_k^\top Q_k x_k + u_k^\top R_k u_k\right),
\end{multline}
where $Q_k\succeq 0$ and $R_k\succ 0$ are the weight matrices.
For each $t$, let $\Pi_t$ denote the class of deterministic causal policies $\pi=(\pi_t,\ldots,\pi_{N-1})$, where $\pi_k:(\bbR^n)^{k+1}\times (\bbR^m)^{k} \to \bbR^m$ is Borel measurable for $k=t,\ldots,N-1$.
Under $\pi\in\Pi_t$, the control input is generated as 
$u_k=\pi_k(h_k)$ before the realization of $w_k$, where $h_k=(x_{0:k},u_{0:k})$ is the information available to the controller.
The finite-horizon risk-sensitive linear quadratic regulator (LQR) problem is formulated as 
\begin{align}\label{eq:LEQRProb}
    J^* = \inf_{ \pi\in\Pi_0} \frac{1}{\theta}\log \bbE \left[e^{\theta L_0(x_{0:N},u_{0:N-1})} \right],
\end{align}
where $\theta\in\bbR\setminus \{0\}$ is a risk-sensitivity parameter.
\end{subequations}

The parameter $\theta$ determines the degree of risk-sensitivity. 
For example, when $\theta>0$, the exponential transformation of the quadratic cost functional in \eqref{eq:LEQRProb} assigns heavier weight to low-probability high-cost trajectories, resulting in a risk-averse behavior that is expected to reduce the worst-case cost. 
Conversely, setting $\theta<0$ assigns lighter weight on such trajectories, reflecting a more aggressive risk-seeking behavior that may perform better on average. 
Lastly, as $\theta \to 0$, risk-sensitive LQR recovers the standard risk-neutral LQR problem, implying that it generalizes LQR.

From \eqref{eq:LEQRProb}, we define the optimal cost-to-go function at each time $t=0,\ldots,N-1$ and state $x\in\bbR^n$ by
\begin{align}\label{eq:costToGo}
    V_t (x) \coloneqq \inf_{\pi\in\Pi_t} \frac{1}{\theta} \log \bbE \left[ e^{\theta L_t(x_{t:N},u_{t:N-1})} \mid x_t = x \right],
\end{align}
where we let $V_N(x) \coloneqq x^\top Q_N x$ for notational consistency.
In particular, $J^* = V_0(x^\ini)$.
As established in \cite{Jaco73,Whit81}, the optimal policy of the risk-sensitive LQR problem admits a linear state feedback form and the cost-to-go function takes a quadratic form, as summarized in the following lemma.
\begin{lem1}\upshape\label{lem:LEQR}
    For a given $\theta\in\bbR\setminus \{0\}$, suppose that the backward Riccati recursion 
    \begin{subequations}\label{eq:riccati}
    \begin{align}
        \!\!P_N \!&=\! Q_N, ~~ \eta_N =0, \label{eq:riccatiIni} \\
        \!\!P_t \!&=\! Q_t + A_t^\top\tilde{P}_{t+1}A_t \!-\! A_t^\top \tilde{P}_{t+1}B_t\tilde{H}_t^{-1}B_t^\top \tilde{P}_{t+1}A_t , \label{eq:riccatiP}\\
        \!\!\eta_t \!&=\! \eta_{t+1} - \frac{1}{2\theta} \log |S_t W_t|,  \label{eq:riccatiEta}\\
        \!\!\tilde{P}_{t+1} \!&=\! P_{t+1}+2\theta P_{t+1}S_t^{-1}P_{t+1} , \label{eq:riccatiPtil}\\
        \!\!\tilde{H}_t \!&=\! R_t  + B_t^\top \tilde{P}_{t+1}B_t,\label{eq:riccatiHtil}\\
        \!\!S_t \!&=\! W_t^{-1}-2\theta P_{t+1}\label{eq:riccatiS},
    \end{align}    
    \end{subequations}
    is well-defined for $t=0,\ldots, N-1$ and satisfies
    \begin{align}\label{eq:feasibility}
        S_t \succ 0, ~~~~ \forall t=0,\ldots, N-1.
    \end{align}
    Then, the optimal policy $\pi^*$ of the finite-horizon risk-sensitive LQR problem \eqref{eq:LEQR} takes the state feedback form
    \begin{align}
        u_t^* &= \pi_t^*(h_t) = \tilde{K}_t x_t, \nonumber \\
        \tilde{K}_t &\coloneqq - \left(R_t + B_t^\top\tilde{P}_{t+1}B_t \right)^{-1} \!B_t^\top \tilde{P}_{t+1}A_t. \label{eq:LEQRinputgain}
    \end{align}
    Moreover, $P_t\succeq 0$ for all $t=0,\ldots, N-1$, and the optimal cost-to-go function is given by $V_t(x) = x^\top P_t x + \eta_t$
    for all $x\in\bbR^n$ and $t=0,\ldots,N$.
\end{lem1}
\begin{proof}
    See Appendix~\ref{apdx:LEQR}.
\end{proof}

Lemma~\ref{lem:LEQR} shows that certainty equivalence does not hold for risk-sensitive LQR.
Certainty equivalence refers to the property that the optimal feedback gain is unaffected by the statistics of the process noise.
This property allows the gain to be computed for the noise-free system and applied to the noisy one without any loss of optimality.
In risk-sensitive LQR, however, the process noise covariance $W_t$ explicitly appears in the Riccati recursion \eqref{eq:riccati}
through $S_t$ and $\tilde{P}_{t+1}$, and hence, the optimal gain $\tilde{K}_t$ needs to be recomputed whenever the process noise statistics change.

The feasibility condition \eqref{eq:feasibility} essentially bounds the intensity of the process noise and the degree of risk-sensitivity.
This is because the condition is more likely to fail as either $\theta$ or the spectral radius of $W_t$ increases.
This is rather natural since the integral of an exponential quadratic functional in \eqref{eq:LEQRProb} may diverge when the process noise is too strong or high-cost realizations are weighted too heavily.
Consequently, the failure of this condition indicates the loss of feasibility of the optimization problem rather than the loss of gain invariance.
In fact, the condition is always satisfied for $\theta<0$, implying that the optimal gain is always well-defined in the risk-seeking case.

\section{Main Results}\label{sec:main}

Input randomization, which refers to deliberately adding random noise into the nominal control input, has been employed in various applications for ensuring privacy or promoting exploration.
Yet, its effects on the optimal policy and cost in risk-sensitive LQR remain unexplored.
In what follows, we formulate the input-randomized risk-sensitive LQR problem and show that the \textit{optimal gain remains invariant to the input noise covariance}.
This obviates the need to redesign the controller even when the input noise covariance changes, potentially simplifying application specific optimization problems.

\subsection{Input-randomized risk-sensitive LQR problem}

To model input randomization, we decompose $u_t$ as
\begin{align}\label{eq:inputDecomp}
    u_t = \bar{u}_t + e_t,
\end{align}
where $\bar{u}_{t}$ denotes the nominal control input to be optimized, and $e_t\in\bbR^m$ is the input noise injected by the controller.
We assume that $e_t$ is mutually independent zero-mean Gaussian with covariance $\Sigma_t \succ 0$, i.e., $e_t\sim \calN(0,\Sigma_t)$, and also independent of the process noise $w_t$.
Importantly, $\bar{u}_t$ is determined before $e_t$ is realized, and hence, cannot directly compensate for it.

The dynamics \eqref{eq:sys} under \eqref{eq:inputDecomp} can be rewritten as
\begin{subequations}\label{eq:LEQR2}
\begin{align}\label{eq:sys2}
    x_{t+1} = A_tx_t + B_t\bar{u}_t + (w_t + B_te_t),
\end{align}
where $\bar{w}_t \coloneqq w_t + B_te_t$ can be interpreted as the effective process noise.
Accordingly, the input-randomized risk-sensitive LQR problem is formulated as
\begin{align}\label{eq:LEQRProb2}
    \bar{J}^* = \inf_{ \bar{\pi} \in \Pi_0} \frac{1}{\theta}\log \bbE \left[e^{\theta L_0(x_{0:N},\bar{u}_{0:N-1}+e_{0:N-1})} \right],
\end{align}
\end{subequations}
where $\bar{u}_t = \bar{\pi}_t(h_t)$ under a policy $\bar{\pi}=(\bar{\pi}_0,\ldots,\bar{\pi}_{N-1})\in\Pi_0$.
Note that the expectation is taken with respect to both $w_{0:N-1}$ and $e_{0:N-1}$, and that the cost penalizes the realized noisy input $u_t=\bar{u}_t+e_t$, not just the nominal input $\bar{u}_t$.
As in \eqref{eq:costToGo}, we define the optimal cost-to-go function at each time $t=0,\ldots,N-1$ and state $x\in\bbR^n$ by
\begin{align}\label{eq:costToGo2}
    \!\!\bar{V}_t (x) \!\coloneqq\! \inf_{ \bar{\pi} \in \Pi_t} \! \frac{1}{\theta} \log \bbE \! \left[ e^{\theta L_t(x_{t:N},\bar{u}_{t:N-1}+e_{t:N-1})} \!\mid\! x_t \!=\! x \right]\!,
\end{align}
where we let $\bar{V}_N(x) \coloneqq x^\top Q_N x$.

\begin{rem1}\upshape\label{rem:leo}

The effect of input randomization in risk-sensitive linear quadratic Gaussian (LQG) control has been studied in \cite{LleoRung25}, where input noise is introduced to promote exploration.
Our formulation differs in how the input cost is evaluated.
Specifically, \cite{LleoRung25} averages the quadratic cost over the input noise distribution before it enters the exponential function. 
In our notation, this replaces $(\bar{u}_k+e_k)^\top R_k (\bar{u}_k+e_k)$ associated with \eqref{eq:costQuad} and \eqref{eq:LEQRProb2} by
\begin{align*}
    \bbE\left[ (\bar{u}_k+e_k)^\top R_k (\bar{u}_k+e_k) \right] = 
    \bar{u}_k^\top R_k \bar{u}_k + \Tr(R_k\Sigma_k),
\end{align*}
where $\Tr(\cdot)$ denotes the trace operation. 
The trace term is independent of the policy. 
However, input noise still enters the state dynamics, and thus, alters the optimal gain as shown in Lemma~\ref{lem:LEQR}.
In our formulation, the exponential function penalizes the realized input $\bar{u}_t+e_t$ directly, thereby accounting for the actual control effort.
This structure will be shown to attribute to the invariance of the optimal gain.
\end{rem1}

\subsection{Motivating example: Dynamic watermarking for replay attack detection}\label{sec:watermark}

We provide a motivating example in which input randomization arises naturally and gain invariance simplifies the associated optimization problem.
Suppose that the state measurements of the plant \eqref{eq:sys} are vulnerable to a replay attack.
Specifically, let $z_t$ denote the state packet received by the controller at time $t$, with $z_t=x_t$ during normal operation.
During an attack interval $[t_a^-,t_a^+]$, an adversary replays previously recorded measurements, as $ z_t=x_{t-\tau}$, where $\tau>0$ is the replay delay.
Since the replayed states were generated by the plant, they may appear statistically plausible and be difficult to distinguish from the true current measurements.

Dynamic watermarking is a representative countermeasure against replay attacks \cite{MoyiSino09,SatcKuma17,NahaTeix23}.
It injects a private excitation signal, called watermark, into the plant input, whose realization is known to the controller but remains hidden from the adversary.
In our setting, the input noise $e_t$ in \eqref{eq:inputDecomp} corresponds exactly to this watermark. 

To see how the watermark enables attack detection, we first consider the time-invariant version of system \eqref{eq:sys} obtained by omitting the subscript $t$ for simplicity, and define the residual as $r_{t+1} \coloneqq z_{t+1} - Az_t -B( \bar{u}_t + e_t )$.
Under normal operation, $r_{t+1}=w_t$, and is therefore uncorrelated with $e_t$. 
When both $z_t$ and $z_{t+1}$ are replayed, i.e., $t,t+1\in [t_a^-, t_a^+]$, 
\begin{align*}
    r_{t+1} &= 
        w_{t-\tau} + B\left(\bar{u}_{t-\tau} +e_{t-\tau} -\bar{u}_t-e_t \right).
\end{align*}
Since $e_t$ is independently sampled and $\bar u_t$ is chosen before $e_t$ is realized, it follows that $\bbE[r_{t+1} e_t^\top] = -B\Sigma$.
Consequently, a replay attack induces a nonzero correlation between $r_{t+1}$ and $e_t$, which can be detected, for example, using a $\chi^2$ detector.

Increasing $\Sigma$ generally enlarges $B\Sigma$ and makes detection easier, but it also degrades control performance.
This motivates choosing $\Sigma$ to maximize the detection probability subject to a bound on the resulting performance loss:
\begin{align}\label{eq:optProb}
    \sup_{\Sigma}\quad & \beta(\Sigma)~~
    \text{s.t.}~~ \bar{J}^{*}-J^{*}\leq\delta,
\end{align}
where $\delta>0$ is a design parameter and $\beta(\Sigma)$ denotes the detection probability determined by the chosen detector and false alarm rate.

In standard LQR/LQG settings, certainty equivalence ensures that the optimal gain remains invariant to the optimization variable $\Sigma$. 
This allows $\bar{J}^*$ to be evaluated under the same controller, simplifying the evaluation of the performance loss $\bar{J}^{*}-J^{*}$. 
Without this property, each candidate $\Sigma$ requires recomputing the optimal gain before $\bar{J}^*$ can be evaluated, turning \eqref{eq:optProb} into a nested optimization problem.
This may have been the reason why related literature evolved mostly around LQR/LQG settings \cite{MoyiSino09,SatcKuma17,NahaTeix23}.
\textit{This motivates our study on input-randomized risk-sensitive LQR, as a first step toward bridging the gap between risk-sensitive control and dynamic watermarking.}

\subsection{Invariance of the optimal gain}
Since $\bar{w}_t \sim \calN(0,W_t+B_t \Sigma_t B_t^\top)$, one might expect the optimal gain to change with $\Sigma_t$ because, as discussed after Lemma~\ref{lem:LEQR}, certainty equivalence does not hold for risk-sensitive LQR.
The following theorem, however, shows that the optimal gain is invariant under input randomization, while the optimal cost-to-go function changes only by a constant that depends on $\Sigma_t$.

\begin{thm1}\upshape\label{thm:invariance}
    For a given $\theta\in\bbR\setminus \{0\}$, suppose that the backward Riccati recursion in \eqref{eq:riccati} is well-defined for $t=0,\ldots, N-1$ and satisfies \eqref{eq:feasibility}. 
    In addition, assume that 
    \begin{align}\label{eq:feasibility2}
        \Xi_t\coloneqq\Sigma_t^{-1} - 2\theta \tilde{H}_t \succ 0, ~~~~ \forall t=0,\ldots, N-1.
    \end{align}
    Then, the following statements hold:
    \begin{enumerate}[leftmargin=*]
        \item The optimal nominal policy $\bar{\pi}^*$ of the finite-horizon input-randomized risk-sensitive LQR problem \eqref{eq:LEQR2} takes the state feedback form $\bar{u}_t^* = \bar{\pi}_t^*(h_t) = \tilde{K}_t x_t$,
        where $\tilde{K}_t$ is defined in \eqref{eq:LEQRinputgain}.
        \item The optimal cost-to-go function is given by $\bar{V}_t(x) = x^\top P_t x + \bar{\eta}_t$ for all $t=0,\ldots, N$, where $\bar{\eta}_N = 0$ and 
        \begin{align*}
            \bar{\eta}_t \!&=\! \bar{\eta}_{t+1} \!-\! \frac{1}{2\theta}\left( \log (|S_tW_t| \!\cdot\! |\Xi_t\Sigma_t|) \right), ~ t=0,\ldots, N-1. 
        \end{align*}
    \end{enumerate}
\end{thm1}

\begin{proof}
    For notational clarity, within this proof, we write $\bbE_{w,e}[\cdot]$, $\bbE_{w}[\cdot]$, and $\bbE_{e}[\cdot]$ to denote expectations with respect to the corresponding noise sequences $(w_{t:N-1},e_{t:N-1})$, $w_{t:N-1}$, and $e_{t:N-1}$, respectively.

    We prove by backward induction that $\bar{V}_t(x) = x^\top P_t x + \bar{\eta}_t$.
    Since $\bar{V}_N(x) = x^\top Q_N x$, the claim holds at $t=N$ with $P_N=Q_N\succeq 0$ and $\bar{\eta}_N=0$. 
    Suppose now that the claim holds at $t+1$ for some $t=0,\ldots, N-1$.
    By Bellman's principle of optimality, 
    \begin{align} \label{eq:barPsi}
       \bar{V}_t(x) &= \inf_{\bar{u}} \frac{1}{\theta}\log \bbE_{w,e}\Big[ e^{\theta\left(x^\top Q_tx + \left(\bar{u}+e_t \right)^\top R_t \left(\bar{u}+e_t \right)\right)} \\
       & \qquad\qquad\qquad\qquad \quad\times e^{\theta\bar{V}_{t+1}(A_tx + B_t(\bar{u}+e_t )+w_t)} \Big ] \nonumber\\
       &=\bar{\eta}_{t+1} + \inf_{\bar{u}} \frac{1}{\theta} \log \bbE_{e}\Big[e^{\theta\left(x^\top Q_tx + \left(\bar{u}+e_t \right)^\top R_t \left(\bar{u}+e_t \right)\right)} \nonumber \\
       & \!\!\!\! \!\!\!\! \!\!\!\!\times \bbE_w\left[ e^{\theta (A_tx+B_t(\bar{u}+e_t)+w_t)^\top P_{t+1} (A_tx+B_t(\bar{u}+e_t)+w_t)}\right] \Big ]. \nonumber
    \end{align} 
    For fixed $\bar{u}$ and $e_t$, Lemma~\ref{lem:GaussInt} is applicable to the inner expectation of \eqref{eq:barPsi} under \eqref{eq:feasibility}, leading to
    \begin{align}\label{eq:barPsi2}
        &\bar{V}_t(x) = \bar{\eta}_{t+1} - \frac{1}{2\theta} \log(|S_tW_t|) + \inf_{\bar{u}} \frac{1}{\theta} \log \bbE_{e}\Big[e^{\theta x^\top Q_tx } \nonumber \\
        & \!\times e^{\theta \left( \left(\bar{u}+e_t \right)^\top R_t \left(\bar{u}+e_t \right) + (A_tx + B_t(\bar{u}+e_t) )^\top \tilde{P}_{t+1} (A_tx + B_t(\bar{u}+e_t) ) \right)}\Big ] \nonumber \\
        &\!=\bar{\eta}_{t+1} \!-\! \frac{1}{2\theta} \!\log(|S_tW_t|) \!+\! \inf_{\bar{u}} \frac{1}{\theta}\! \log \bbE_{e}\!\Big[e^{\theta x^\top(Q_t + A_t^\top \tilde{P}_{t+1}A_t) x} \nonumber\\
        & \qquad \qquad \qquad \times e^{\theta ((\bar{u}+e_t)^\top \tilde{H}_t (\bar{u}+e_t) +2(\bar{u}+e_t)^\top B_t^\top \tilde{P}_{t+1}A_tx ) } \Big] \nonumber \\
        &\!= x^\top P_t x + \bar{\eta}_{t+1}- \frac{1}{2\theta} \log(|S_tW_t|) \\
        & \qquad \qquad \quad  + \inf_{\bar{u}} \frac{1}{\theta} \log \bbE_e \Big[ e^{\theta (\bar{u}-\tilde{K}_tx + e_t)^\top \tilde{H}_t(\bar{u}-\tilde{K}_tx + e_t)  } \Big]. \nonumber
    \end{align}
    By Lemma~\ref{lem:LEQR}, $\tilde{P}_{t+1}\succeq 0$, and hence, $\tilde{H}_t \succ 0$ by \eqref{eq:riccatiHtil} and $R_t\succ 0$.
    Therefore, under \eqref{eq:feasibility2}, Lemma~\ref{lem:GaussInt} is again applicable to the expectation in the last equation of \eqref{eq:barPsi2}, yielding
    \begin{align}\label{eq:barPsi3}
        \bar{V}_t(x) = x^\top P_t x + \bar{\eta}_{t+1} &- \frac{1}{2\theta}\left( \log (|S_tW_t| \cdot |\Xi_t\Sigma_t|) \right) \\
        & + \inf_{\bar{u}}(\bar{u}-\tilde{K}_tx)^\top \tilde{G}_t(\bar{u}-\tilde{K}_tx), \nonumber
    \end{align}
    where $\tilde{G}_t\coloneqq \tilde{H}_t + 2\theta \tilde{H}_t \Xi_t^{-1} \tilde{H}_t$.
    Indeed, $\tilde{G}_t\succ 0$ because $\tilde{H}_t\succ 0$ and $\Xi_t\succ 0$.
    Consequently, the minimizer is uniquely attained at $\bar{u}^* = \tilde{K}_t x$, and thus, $\bar{V}_t(x) = x^\top P_t x + \bar{\eta}_{t+1} - \frac{1}{2\theta} \log (|S_tW_t| \cdot |\Xi_t \Sigma_t|).$
    This concludes the proof.
\end{proof}

Theorem~\ref{thm:invariance} shows that the input noise covariance $\Sigma_t$ neither enters the Riccati recursion \eqref{eq:riccati} nor alters the optimal gain.
This is far from obvious, as it inflates the effective process noise covariance, which generally affects the optimal gain in risk-sensitive LQR, as shown in Lemma~\ref{lem:LEQR}.
The invariance is attributed to the fact that the input noise $e_t$ enters not only the system dynamics but also the quadratic input cost in \eqref{eq:LEQRProb2}.
As seen in \eqref{eq:barPsi3}, these two effects together reshape the curvature of the exponential quadratic function from $\tilde{H}_t$ to $\tilde{G}_t$, while leaving the minimizer at $\bar u_t^*=\tilde K_t x$ unchanged.

\subsection{Cost increment and joint feasibility}
Although the optimal gain remains invariant, input randomization changes the optimal cost through the additive term $\bar{\eta}_t$, whose recursion additionally accumulates $-\frac{1}{2\theta}\log|\Xi_t\Sigma_t|$ relative to \eqref{eq:riccatiEta}.
Nonetheless, because the gain is invariant, the performance loss $\bar{J}^* - J^*$ can be computed readily, without recomputing the Riccati recursion.
This could be particularly useful when performance loss has to be computed across varying values of $\Sigma_t$, as in \eqref{eq:optProb}.


\begin{prop1}\upshape 
    For a given $\theta \in \bbR \setminus \{0\}$, suppose that the assumptions of Theorem~\ref{thm:invariance} hold.
    Then, for any initial value $x^\ini$, the performance loss due to input randomization is  
    \begin{align*}
        \bar{J}^* - J^* = -\frac{1}{2\theta} \textstyle\sum_{t=0}^{N-1} \log |\Xi_t \Sigma_t| >0.
    \end{align*}
\end{prop1}
\begin{proof}
    By definition, $J^* = V_0(x^\ini)$ and $\bar{J}^* = \bar{V}_0(x^\ini)$ (see \eqref{eq:LEQRProb}, \eqref{eq:costToGo}, \eqref{eq:LEQRProb2}, and \eqref{eq:costToGo2}). 
    Therefore, by Lemma~\ref{lem:LEQR} and Theorem~\ref{thm:invariance}, it is obtained that 
    \begin{align}\label{eq:costDiff}
        \bar{J}^* - J^*
        &= \bar{\eta}_0 - \eta_0 = -\frac{1}{2\theta}\textstyle\sum_{t=0}^{N-1}  \log |\Xi_t \Sigma_t|.
    \end{align}
    Since $\Xi_t\succ 0$ by \eqref{eq:feasibility2} and $\Sigma_t \succ 0$, we have $\Sigma_t^{1/2}\Xi_t \Sigma_t^{1/2} \succ 0$.
    Therefore,
    $|\Xi_t\Sigma_t| = |\Sigma_t^{1/2}\Xi_t \Sigma_t^{1/2}|>0$, and hence, the logarithms in \eqref{eq:costDiff} are well-defined. 
    
    Define $M_t\coloneqq\Sigma_t^{1/2} \tilde{H}_t \Sigma_t^{1/2}$, so that $\Sigma_t^{1/2}\Xi_t \Sigma_t^{1/2} = I-2\theta M_t$.
    As determinant equals the product of all eigenvalues, 
    \begin{align}\label{eq:detProd}
        |\Xi_t\Sigma_t| = \textstyle \prod_{i=1}^m (1 - 2\theta \lambda_{i,t})>0,
    \end{align}
    where $\lambda_{1,t},\ldots,\lambda_{m,t}$ denote the eigenvalues of $M_t$.
    Since $M_t$ is congruent to $\tilde{H}_t\succ 0$, we have $\lambda_{i,t}>0$ for all $i=1,\ldots, m$.
    Thus, $1-2\theta \lambda_{i,t} >1$ for $\theta<0$ and $1-2\theta \lambda_{i,t}<1$ for $\theta>0$.
    Meanwhile, positive definiteness of
    $I-2\theta M_t$ additionally gives $1-2\theta \lambda_{i,t} >0$ regardless of the sign of $\theta$.
    Therefore, it follows from \eqref{eq:detProd} that, for any $\theta\in\bbR \setminus \{0\}$, $\frac{1}{2\theta}\log|\Xi_t\Sigma_t|<0$ for all $t=0,\ldots, N-1$, leading to $\bar{J}^*-J^*>0$. 
    This concludes the proof.
\end{proof}

It should be noted that the additional feasibility condition \eqref{eq:feasibility2} arises due to input randomization.
Similar to the condition \eqref{eq:feasibility}, it restricts the admissible intensity of the input noise $e_t$, and can only fail in the risk-averse regime ($\theta>0$).
Characterizing the joint feasibility of \eqref{eq:feasibility} and \eqref{eq:feasibility2} may appear difficult, as $\Xi_t$ depends on $\theta$ both directly and indirectly through the Riccati recursion \eqref{eq:riccati}.
The following lemma, however, reveals that both conditions are feasible when $\theta$ is chosen sufficiently small.

\begin{lem1}\upshape\label{lem:monotone}
    Define 
    \begin{align*}
        \Theta \!\coloneqq\! \left\{\theta\in\bbR\setminus\{0\} \mid\eqref{eq:riccati}~\mbox{is well-defined},~\eqref{eq:feasibility}~\mbox{and}~\eqref{eq:feasibility2}~\mbox{hold} \right\}.
    \end{align*}
Then, there exists a constant $\bar{\theta}\in (0,\infty)$ such that $\Theta = (-\infty, \bar{\theta}) \setminus \{0\}$.
\end{lem1}

\begin{proof}
Throughout the proof, the dependence on $\theta$ is made explicit for the quantities generated by \eqref{eq:riccati}, for example, by writing $P_t(\theta)$.
It suffices to show that $\Theta$ (i) contains a strictly positive element; (ii) is monotone, i.e., $\theta_2\in\Theta$ implies $\theta_1\in\Theta$ for all nonzero $\theta_1< \theta_2$; (iii) is bounded above. 
Since $\Theta$ is open, as \eqref{eq:feasibility} and \eqref{eq:feasibility2} vary continuously in $\theta$, these properties give $\bar{\theta}=\sup \Theta \in (0,\infty)$.

The positive definiteness of $(W_t,\Sigma_t)$ and positive semi definiteness of $(P_{t+1}, \tilde{H}_t )$ alone do not ensure joint feasibility of \eqref{eq:feasibility} and \eqref{eq:feasibility2} for some positive $\theta$, as the latter depend on $\theta$. 
Therefore, to establish this, we extend the Riccati recursion \eqref{eq:riccati} continuously to $\theta = 0$, where it reduces to the standard LQR Riccati recursion with
\begin{align*}
    S_t(0) = W_t^{-1} \succ 0, ~~ \Xi_t(0)=\Sigma_t^{-1} \succ 0, ~~ t=0,\ldots,N-1.
\end{align*}
Since the horizon is finite, it follows from the continuity that there exists a uniform $\epsilon>0$ such that \eqref{eq:feasibility} and \eqref{eq:feasibility2} hold for all $\theta \in (0,\epsilon)$.

To show monotonicity, let $0 < \theta_1 < \theta_2$, and suppose that $\theta_2\in\Theta$.
We only consider positive $\theta_1$ because $(-\infty, 0) \subset \Theta$ follows immediately from the positive definiteness of $(W_t,\Sigma_t)$ and positive semi definiteness of $(P_{t+1}, \tilde{H}_t )$.
We show by backward induction that $\theta_1\in\Theta$ and $0\preceq \theta_1 P_t(\theta_1) \preceq  \theta_2 P_t(\theta_2)$ for all
    $t = 0,\ldots,N$.
    The claim holds at $t = N$ since $P_N(\theta_1)=P_N(\theta_2) = Q_N\succeq 0$. 
    Now, suppose that $0 \preceq \theta_1 P_{t+1}(\theta_1) \preceq \theta_2 P_{t+1}(\theta_2)$ for some $t=0,\ldots,N-1$. 
    Then, 
    \begin{align*}
        S_t(\theta_1) &= W_t^{-1} - 2\theta_1 P_{t+1}(\theta_1) \\
        &\succeq W_t^{-1} - 2\theta_2 P_{t+1}(\theta_2) = S_t(\theta_2) \succ 0,
    \end{align*}
    and thus, \eqref{eq:feasibility} holds at $t$ for $\theta_1$.
Moreover, since $S_t(\theta_1)\succ 0$, Lemma~\ref{lem:GaussInt} gives $\tilde{P}_{t+1}(\theta_1)\succeq 0$. 
As a result, $\tilde{H}_t(\theta_1)\succ 0$, and the recursion \eqref{eq:riccati} is well-defined at $t$ for $\theta_1$.

\begin{figure*}[t]
    \centering
        \input{fig/plot}
    \vspace{-0.3cm}
    \caption{ Risk-sensitive cost $\bar{J}(K;\Sigma)$ in \eqref{eq:rscost} versus the gain $K$ for three choices of the input cost term $\rho(K,e_0;\Sigma)$ and three values of the input noise covariance $\Sigma$. 
    The filled circles denote the minimum of each plot.}
    \label{fig:scalar}
\end{figure*}

We now turn our attention to verifying \eqref{eq:feasibility2}. 
Since the Loewner order is not preserved under matrix multiplication \cite[Section 7.7]{HornJohn12}, i.e., $X_1\succeq X_2$ and $Y_1\succeq Y_2$ do not imply $X_1Y_1\succeq X_2Y_2$, we use the following alternative expression for $2\theta\tilde{P}_{t+1}(\theta)$: 
    \begin{align*}
        &2\theta \tilde{P}_{t+1}(\theta) = 2\theta P_{t+1}(\theta) + 4\theta^2 P_{t+1}(\theta) S_t^{-1}(\theta) P_{t+1}(\theta) \nonumber \\
        &= (W_t^{-1} - S_t(\theta)) + (W_t^{-1} - S_t(\theta)) S_t^{-1}(\theta) (W_t^{-1} \!-\! S_t(\theta)) \nonumber \\
        &=-W_t^{-1} + W_t^{-1} S_t^{-1}(\theta) W_t^{-1}.
    \end{align*}
Because inversion reverses the Loewner order and congruence preserves it \cite[Theorem 7.7.2 and Corollary 7.7.4]{HornJohn12}, it is obtained that $W_t^{-1} S_t^{-1}(\theta_1) W_t^{-1} \preceq W_t^{-1} S_t^{-1}(\theta_2) W_t^{-1}$, and thus, $0 \preceq 2\theta_1 \tilde{P}_{t+1}(\theta_1) \preceq 2\theta_2 \tilde{P}_{t+1}(\theta_2)$.
    Consequently, by \eqref{eq:riccatiHtil} and $R_t\succ 0$, we have
    \begin{align*}
        \Xi_t(\theta_1) &= \Sigma_t^{-1} - 2\theta_1(R_t  + B_t^\top \tilde{P}_{t+1}(\theta_1)B_t) \\
        &\succ \Sigma_t^{-1} - 2\theta_2(R_t  + B_t^\top \tilde{P}_{t+1}(\theta_2)B_t) = \Xi_t(\theta_2) \succ 0,
    \end{align*}
    and \eqref{eq:feasibility2}  holds at $t$ for $\theta_1$ as well.

Let $\Gamma_t(\theta)$ denote the block matrix in \eqref{eq:Gamma}, with its dependence on $\theta$ made explicit, so that $ \theta P_{t}(\theta)$ is the Schur complement of the block $\theta \tilde{H}_t(\theta)$ in $\theta \Gamma_t(\theta)$.
Since $\theta_1 < \theta_2$ and $\theta_1\tilde{P}_{t+1}(\theta_1)\preceq \theta_2\tilde{P}_{t+1}(\theta_2)$, it follows that $\theta_1 \Gamma_t(\theta_1)\preceq \theta_2 \Gamma_t(\theta_2)$.
As the Schur complement is monotone with respect to the Loewner order \cite[Corollary~1]{AndeTrap75}, $\theta_1 P_t(\theta_1)\preceq  \theta_2 P_t(\theta_2)$.
In addition, $0 \preceq \theta_1 P_t(\theta_1)$ by Lemma~\ref{lem:LEQR}. 
Together, these establish that $\theta_1\in\Theta$ by backward induction.

Finally, we show that $\Theta$ is bounded above. 
By definition, $\Xi_t(\theta) = \Sigma_t^{-1} - 2\theta \tilde{H}_t(\theta) \preceq \Sigma_t^{-1} -2\theta R_t$.
Since $\Sigma_t$ and $R_t$ are independent of $\theta$, the right-hand-side fails to be positive definite when $\theta > \lambda_{\min}(\Sigma_t^{-1}) / 2\lambda_{\min}(R_t)$,
where $\lambda_{\min}(\cdot)$ denotes the minimum eigenvalue. 
Therefore, $\Theta$ is bounded above, and this concludes the proof.
\end{proof}

Lemma~\ref{lem:monotone} reveals that the joint feasibility of \eqref{eq:feasibility} and \eqref{eq:feasibility2} is monotone in $\theta$.
Thus, if the conditions hold for a given $\theta$, they remain valid for all smaller values.
Moreover, verifying the conditions for a fixed $\theta$ is direct once the recursion \eqref{eq:riccati} has been computed.

\subsection{Numerical simulation}\label{sec:simul}

Consider a one-step ($N=1$) scalar case ($n=m=1$) of the input-randomized system \eqref{eq:sys2}, given by
\begin{align*}
    x_1=x_0+(\bar{u}_0+e_0) + w_0,    
\end{align*}
where $x_0=1$, $w_0\sim\calN(0,0.1)$, and
$e_0\sim\calN(0,\Sigma)$.
Since $x_0$ is fixed and $N=1$, every policy $\pi\in\Pi_0$ specifies a scalar nominal input $\bar{u}_0=\pi
(x_0)=Kx_0$ for some $K\in\bbR$.
Accordingly, the risk-sensitive cost is defined as 
\begin{align}\label{eq:rscost}
    \bar{J}(K;\Sigma)
\coloneqq  \frac{1}{\theta}\log\bbE\!\left[
e^{\theta\left(Q_1x_1^2+Q_0x_0^2+R_0\rho(K,e_0;\Sigma)\right)}
\right],
\end{align}
where we set $\theta=0.2$, $Q_0=0$, and $Q_1=R_0=1$.
We compare three choices of the input cost term:
\begin{align*}
    \rho(K,e_0;\Sigma) = \begin{cases}
        (K+e_0)^2, & \mbox{(Case 1)} \\
        \bbE\left[(K+e_0)^2\right] = K^2 + \Sigma, & \mbox{(Case 2)} \\
        K^2. & \mbox{(Case 3)}
    \end{cases}
\end{align*}
Specifically, we examine how the cost $\bar{J}(K;\Sigma)$ depends on the gain $K$ and input noise covariance $\Sigma$.

Case 1 penalizes the realized input $\bar{u}_0+e_0 = K+e_0$ directly, as in our proposed formulation \eqref{eq:LEQR2}. 
Indeed, the optimal cost in \eqref{eq:LEQRProb2} is obtained as $\bar{J}^*=\inf_{K\in\bbR} \bar{J}(K;\Sigma)$.
Case 2 follows \cite{LleoRung25} by averaging the quadratic input cost over $e_0$ before applying the exponential function, as described in Remark~\ref{rem:leo}. 
Finally, Case 3 only penalizes the nominal control input $\bar{u}_0=K$ and omits $e_0$ from the input cost term.

The cost $\bar{J}(K;\Sigma)$ admits a closed-form expression for each case. 
For Case 1 as an example,
since $N=1$ and $\bar{\eta}_1=0$, evaluating the  expression inside the infimum in \eqref{eq:barPsi3} at $\bar{u}=Kx_0$, $x=x_0=1$, and $t=0$ gives
\begin{align*}
    \bar{J}(K;\Sigma) = P_0 - \frac{1}{2\theta} \left( \log (S_0W_0\Xi_0\Sigma) \right) +  \tilde{G}_0(K - \tilde{K}_0)^2.    
\end{align*}
Computing the Riccati recursion \eqref{eq:riccati} and substituting yields 
\begin{align*}
    \bar J(K;\!\Sigma)
\!=\! \frac{25}{49}
\!-\! \frac{5}{2} \!\log(0.96 \!-\! 0.784\Sigma) \!+\!
\frac{49}{24 \!-\! \frac{98}{5}\Sigma}
\!\left(\! K \!+\! \frac{25}{49}\!\right)^2\!\!,
\end{align*}
which can be done analogously for Cases 2 and 3.

Fig.~\ref{fig:scalar} shows that the optimal gain is invariant in Case 1 but varies with $\Sigma$ in Cases 2 and 3. 
This validates the result of Theorem~\ref{thm:invariance} and highlights the role of penalizing the realized noisy input in preserving the optimal gain.



\subsection{ Interpretation as a zero-sum dynamic game}

In this subsection, we establish a connection between the input-randomized risk-sensitive LQR problem \eqref{eq:LEQR2} and a zero-sum dynamic game that yields the same Riccati recursion and optimal control policy.
This provides an intuitive interpretation of the invariance of the optimal gain under input randomization.

We consider a finite-horizon zero-sum dynamic game between a controller and an auxiliary adversary.
The state dynamics is given by
\begin{align*}
    x_{t+1} = A_tx_t + B_t(\hat{u}_t + d_t) + v_t,
\end{align*}
where the controller chooses $\hat{u}_t\in\bbR^m$, and the adversary chooses $d_t\in\bbR^m$ and $v_t\in\bbR^n$.
Specifically, let $h_t=(x_{0:t}, \hat{u}_{0,t-1}, v_{0,t-1}, d_{0,t-1})$ denote the information available to both players at time $t$.
The controller first chooses a deterministic causal policy $\phi_t$ and generates $\hat{u}_t = \phi_t(h_t)$. 
After observing $\hat{u}_t$, the adversary chooses a deterministic causal policy $\psi_t$ and generates $(v_t,d_t) = \psi_t(h_t,\hat{u}_t)$.
We denote by $\Phi_t$ and $\Psi_t$ the respective classes of admissible policies $\phi=(\phi_t,\ldots,\phi_{N-1})$ and $\psi=(\psi_t,\ldots, \psi_{N-1})$ at time $t$. 

The value of the game is defined by  
\begin{multline}\label{eq:game}
    \hat{J}^{*}
    = \inf_{\phi\in \Phi_0} \sup_{\psi\in\Psi_0}
    L_0(x_{0:N}, \hat{u}_{0:N-1} + d_{0:N-1}) \\
    - \frac{1}{2\theta} \sum_{k=0}^{N-1}\left( v_k^\top W_k^{-1}v_k + d_k^\top \Sigma_k^{-1} d_k \right),
\end{multline}
where we restrict our attention to the case $\theta>0$.
Accordingly, the controller's worst-case cost-to-go function at each time $t=0,\ldots,N-1$ and state $x\in\bbR^n$ is defined by
    \begin{align*}
        \hat{V}_t(x) = \inf_{\phi \in \Phi_t} \sup_{\psi \in \Psi_t} L_t(x_{t:N}, \hat{u}_{t:N-1} + d_{t:N-1}) \\
    - \frac{1}{2\theta} \sum_{k=t}^{N-1}\left( v_k^\top W_k^{-1}v_k + d_k^\top \Sigma_k^{-1} d_k \right),
\end{align*}
where $x_t=x$, and we let $\hat{V}_N(x) = x^\top Q_N x$.

In \eqref{eq:game}, the controller chooses $\phi\in\Phi_0$ to minimize the cost, whereas the adversary chooses $\psi\in\Psi_0$ to maximize it. 
Unlike \eqref{eq:LEQR2}, risk sensitivity is not introduced by exponentiating the quadratic trajectory cost $L_0$. 
Instead, the adversarial disturbances $v_t\in\bbR^n$ and $d_t\in\bbR^m$ are penalized through the quadratic terms weighted by $W_t^{-1}$ and $\Sigma_t^{-1}$, respectively.
As $\theta$ increases, these penalties weaken and allows the adversary to inject stronger disturbances, thereby rendering the controller more risk-averse. 

The following theorem shows that the game admits a saddle point whose
controller policy coincides with that of the input-randomized
risk-sensitive LQR problem.

\begin{thm1}\upshape\label{thm:game}
    For a given $\theta >0$, suppose that the backward Riccati recursion in \eqref{eq:riccati} is well-defined for $t=0,\ldots,N-1$, and that the feasibility conditions \eqref{eq:feasibility} and \eqref{eq:feasibility2} hold.
    Then, the zero-sum dynamic game \eqref{eq:game} admits a saddle point\footnote{That is, changing $\phi^*$ alone cannot decrease the objective in \eqref{eq:game}. Likewise, changing $\psi^*$ alone cannot increase it.}
    $(\phi^*,\psi^*)$,
    where the controller's optimal policy takes the state feedback form
    \begin{subequations}\label{eqs:saddle}
    \begin{align}\label{eq:saddleu}
        \hat{u}_t^* = \phi_t^*(h_t) = \tilde{K}_t x_t,
    \end{align}
    while the adversary's optimal policy $(v_t^*,d_t^*) = \psi_t^*(h_t,\hat{u}_t)$ for any controller action $\hat{u}_t$ is given by 
    \begin{align}
         v_t^* &= 2\theta S_t^{-1}P_{t+1}(A_tx_t +B_t(\hat{u}_t+d_t^*)),\label{eq:saddlev} \\
        d_t^* &= 2\theta \Xi_t^{-1}\tilde{H}_t(\hat{u}_t-\tilde{K}_tx_t).\label{eq:saddled} 
    \end{align}
    \end{subequations}
    In particular, $d_t^* = 0$ whenever $\hat{u}_t = \hat{u}_t^*$.
   Moreover, $\hat{V}_t(x) = x^\top P_t x$ for all $t=0,\ldots,N$.
\end{thm1}

\begin{proof}
    We prove by backward induction. 
    Since $\hat{V}_N(x) = x^\top Q_N x$, the claim holds at $t=N$.
    Suppose now that the claim holds at $t+1$ for some $t=0,\ldots, N-1$. 
    Then, by Bellman's principle of optimality,
    \begin{align*}
        &\hat{V}_t(x) = \inf_{\hat{u}} \sup_{(v,d)} ~ x^\top Q_t x + (\hat{u}+d)^\top R_t (\hat{u}+d) \\
        &~~-\frac{1}{2\theta}(v^\top W_t^{-1}v + d^\top \Sigma_t^{-1} d) + (z+v)^\top P_{t+1}(z+v) \\
        &=\inf_{\hat{u}} \sup_{d} ~ x^\top Q_t x + (\hat{u}+d)^\top R_t (\hat{u}+d)-\frac{1}{2\theta} d^\top \Sigma_t^{-1}d \\
        &~~+\sup_v ~ (z+v)^\top P_{t+1}(z+v) -\frac{1}{2\theta} v^\top W_t^{-1}v,
    \end{align*}
    where $z:=A_tx+B_t(\hat{u}+d)$.
    For fixed $\hat{u}$ and $d$, 
    \begin{align*}
        &\sup_v ~ (z+v)^\top P_{t+1}(z+v) -\frac{1}{2\theta} v^\top W_t^{-1}v  \\
        &= \sup_v -\frac{1}{2\theta} (v-2\theta S_t^{-1}P_{t+1}z)^\top S_t  (v - 2\theta S_t^{-1}P_{t+1}z) \\
        & ~~ + z^\top \tilde{P}_{t+1} z \\
        &= (A_tx+B_t(\hat{u}+d))^\top \tilde{P}_{t+1} (A_tx+B_t(\hat{u}+d)),
    \end{align*}
    where the last equation holds since $S_t\succ 0$ under \eqref{eq:feasibility} and the maximizer is obtained at $v^* = 2\theta S_t^{-1}P_{t+1}z$.
    Substituting this yields 
    \begin{align*} 
        &\hat{V}_t(x) \!=\! \inf_{\hat{u}} \sup_{d} ~ x^\top \!(Q_t \!+\! A_t^\top \tilde{P}_{t+1}A_t)x  \!+\! (\hat{u} \!+\! d)^\top \!\tilde{H}_t (\hat{u} \!+\! d)\\
        &~~ + 2(\hat{u}+d)^\top B_t^\top \tilde{P}_{t+1}A_tx - \frac{1}{2\theta}d^\top \Sigma_t^{-1}d \\
        &=x^\top P_t x +\inf_{\hat{u}} \sup_{d} ~(r+d)^\top \tilde{H}_t (r+d) - \frac{1}{2\theta}d^\top \Sigma_t^{-1}d,
    \end{align*}
    where $r:=\hat{u}-\tilde{K}_t x$.
    By completing the square, it is obtained that 
    \begin{align*}
        &(r+d)^\top \tilde{H}_t (r+d) - \frac{1}{2\theta}d^\top \Sigma_t^{-1}d \\
        &= -\frac{1}{2\theta}(d-2\theta \Xi_t^{-1}\tilde{H}_t r)^\top \Xi_t (d-2\theta \Xi_t^{-1}\tilde{H}_t r) + r^\top \tilde{G}_t r, 
    \end{align*}
    where $\tilde{G}_t:=\tilde{H}_t +2\theta \tilde{H}_t \Xi_t^{-1}\tilde{H}_t$.
    Since $\Xi_t\succ 0$ under \eqref{eq:feasibility2}, the maximization over $d$ has value $r^\top \tilde{G}_t r$ and is uniquely attained at $d^*= 2\theta \Xi_t^{-1}\tilde{H}_t r$. Furthermore, $\tilde{G}_t\succ 0$ because $\tilde{H}_t\succ 0$.
Consequently, 
\begin{align*}
    \hat{V}_t(x) = x^\top P_t x + \inf_{\hat{u}} r^\top \tilde{G}_t r = x^\top P_t x, 
\end{align*}
with the minimum uniquely attained at $r=0$, or equivalently, $\hat{u}=\tilde{K}_t x$. 
    This concludes the proof.
\end{proof}

Theorem~\ref{thm:game} offers an intuitive explanation for the gain invariance. 
At the optimal gain $\hat{u}_t^* = \tilde{K}_t x_t$, the adversary's best response is $d_t^*=0$, i.e., it has no incentive to inject any input noise. 
If the controller deviates from the optimal gain, however, the adversary's best response becomes $d_t^*\neq 0$. 
From a game-theoretic perspective, the controller should not deviate from the optimal gain, as doing so only opens the door for the adversary to exploit.

\section{Concluding Remarks}\label{sec:conclusion}
In this paper, we showed that the optimal gain of risk-sensitive LQR is invariant under input randomization when the cost penalizes the realized noisy input.
By exploiting this invariance, we derived a closed-form expression for the resulting increment in the optimal cost and characterized the feasibility of the problem, which holds whenever the risk-sensitivity parameter $\theta$ is sufficiently small. 
We believe that this result can bridge the gap between risk-sensitive control and various applications that employ input randomization.
Future work includes extending the result to the infinite-horizon setting and analyzing the invariance of the gain in input-randomized risk-sensitive LQG problem. 

\bibliographystyle{IEEEtran}
\bibliography{ref}

\appendices

\section{Technical lemma}\label{apdx:GaussInt}

\begin{lem1}\upshape\label{lem:GaussInt}
    Let $n\in\bbN$, and consider matrices $P,W\in\bbR^{n\times n}$ such that $P\succeq 0$ and $W\succ 0$, and a vector $w\in\bbR^n$ such that $w\sim \calN(0,W)$.
    For any $\theta\in\bbR\setminus\{0\}$ such that $S\coloneqq W^{-1}-2\theta P \succ 0$, it holds for every $z\in\bbR^n$ that
    \begin{align*}
        \bbE\left[ e^{\theta (z+w)^\top P (z+w)} \right] = \frac{e^{\theta z^\top \tilde{P} z}}{\sqrt{|SW|}},
    \end{align*}
    where $\tilde{P}\coloneqq P+2\theta PS^{-1}P\succeq 0$.
\end{lem1}

\begin{proof}
    See \cite[Appendix~1]{Jang26}.
\end{proof}

\section{Proof of Lemma~\ref{lem:LEQR}}\label{apdx:LEQR}

We show by backward induction that $P_t\succeq 0$ and $V_t(x) = x^\top P_t x + \eta_t$ for all $t=0,\ldots, N$.
Since $V_N(x) = x^\top Q_N x$, the claim holds at $t=N$ with $P_N=Q_N\succeq 0$ and $\eta_N = 0$.
Suppose now that the claim holds at $t+1$ for some $t=0,\ldots,N-1$. 
Then, by Bellman's principle of optimality and Lemma~\ref{lem:GaussInt}, 
\begin{align}\label{eq:Psi}
    & V_t(x) = \inf_{u} x^\top Q_tx \!+\! u^\top R_t u  \!+\! \frac{1}{\theta}\! \log \bbE\! \left[ e^{\theta V_{t+1}(A_tx+B_tu+w_t)} \right] \nonumber \\
    & = \inf_{u} x^\top Q_tx +u^\top R_t u \nonumber  \\
    &  ~~~~~~~~ +  \frac{1}{\theta} \log \frac{e^{\theta\left(\eta_{t+1}+ (A_tx+B_tu)^\top \tilde{P}_{t+1}(A_tx+B_tu) \right)}}{\sqrt{|S_tW_t|}} \nonumber \\
    & =  x^\top P_t x + \eta_{t} + \inf_{u}  (u - \tilde{K}_tx)^\top \tilde{H}_t (u - \tilde{K}_t x),
\end{align}
where the last equality follows by completing the square in $u$, and is valid for either sign of $\theta$.
Lemma~\ref{lem:GaussInt} also gives $\tilde{P}_{t+1}\succeq 0$, which, together with $R_t\succ 0$, leads to $\tilde{H}_t\succ 0$.
Therefore, the minimizer in \eqref{eq:Psi} is uniquely attained at $u^*=\tilde{K}_t x$ and yields $V_t(x) = x^\top P_t x + \eta_t$, as claimed.
Moreover, consider the block matrix 
\begin{align}\label{eq:Gamma}
    \Gamma_t \coloneqq \begin{bmatrix}
        Q_t + A_t^\top \tilde{P}_{t+1} A_t & A_t^\top \tilde{P}_{t+1} B_t \\ B_t^\top \tilde{P}_{t+1} A_t & \tilde{H}_t
    \end{bmatrix}
\end{align}
which is clearly positive semidefinite, as $Q_t\succeq 0$ and $R_t\succ 0$.
Since $\tilde{H}_t\succ 0$, the Schur complement of $\tilde{H}_t$ in $\Gamma_t$, namely $P_t$, is positive semidefinite. 
This concludes the proof. \hfill $\blacksquare$

\end{document}